\documentclass[a4paper,12pt]{article}
\usepackage{amsmath}
\usepackage{latexsym}
\numberwithin{equation}{section}

\def\R{{\bf R}}

\def\d{\displaystyle}
\def\e{{\varepsilon}}

\def\v#1{\mbox{\boldmath $#1$}}

\newtheorem{thm}{Theorem}[section]

\newtheorem{lem}{Lemma}[section]
\newtheorem{prop}{Proposition}[section]
\newtheorem{rem}{Remark}[section]

\title{Improved Kato's lemma on ordinary differential inequality
and its application to semilinear wave equations}
\author{
Hiroyuki Takamura
\footnote{Department of Complex and Intelligent Systems,
Faculty of Systems Information Science,
Future University Hakodate,
116-2 Kamedanakano-cho,
Hakodate, Hokkaido 041-8655, Japan.
e-mail: takamura@fun.ac.jp.}}
\date{
\[
\begin{array}{ll}
\mbox{\footnotesize{\bf Keywords:}}
& \mbox{\footnotesize ordinary differential inequality, semilinear wave equation, lifespan}\\
\mbox{\footnotesize{\bf MSC2010:}}
& \mbox{\footnotesize primary 35L71, 34A40, 35E15, secondary 35A01, 35B44}\\
\end{array}
\]
}

\begin{document}
\maketitle
\begin{abstract}
We are interested in the upper bound of the lifespan of solutions of semilinear wave equations from above.
For the sub-critical case in high dimensions, it has been believed that the basic tools of its analysis are
Kato's lemma on ordinary differential inequalities and the rescaling argument in the functional method.
But there is a small lack of delicate analysis and no published paper about this.
Here we give a simple alternative proof by means of improved Kato's lemma without any rescaling argument.
\end{abstract}


\section{Introduction}
\par
We consider the initial value problem,
\begin{equation}
\label{IVP}
\left\{
\begin{array}{l}
u_{tt}-\Delta u=|u|^p,\quad \mbox{in}\quad \R^n\times[0,\infty),\\
u(x,0)=\e f(x),\ u_t(x,0)=\e g(x)
\end{array}
\right.
\end{equation}
assuming that $\e>0$ is \lq\lq small."
Let us define a lifespan $T(\e)$ of a solution of (\ref{IVP}) by
\[
T(\e)=\sup\{t>0\ :\ \exists\ \mbox{a solution $u(x,t)$ of (\ref{IVP})
for arbitrarily fixed $(f,g)$.}\},
\]
where \lq\lq solution" means classical one when $p\ge2$.
When $1<p<2$, it means weak one, but sometimes the one
given by associated integral equations to (\ref{IVP})
by standard Strichartz's estimate.
See Sideris \cite{Si84} or Georgiev, Takamura and Zhou \cite{GTZ06} for example on such an argument. 
\par
When $n=1$, we have $T(\e)<\infty$ for any power $p>1$ by Kato \cite{Kato80}.
When $n\ge2$, we have the following Strauss' conjecture on (\ref{IVP})
by Strauss \cite{St81}. 
\[
\begin{array}{lll}
T(\e)=\infty & \mbox{if $p>p_0(n)$ and $\e$ is \lq\lq small"}
& \mbox{(global-in-time existence)},\\
T(\e)<\infty & \mbox{if $1<p\le p_0(n)$}
& \mbox{(blow-up in finite time)},
\end{array}
\]
where $p_0(n)$ is so-called Strauss' exponent defined
by positive root of the quadratic equation, $\gamma(p,n)=0$, where
\begin{equation}
\label{gamma}
\gamma(p,n):=2+(n+1)p-(n-1)p^2.
\end{equation}
That is,
\begin{equation}
\label{p_0(n)}
p_0(n)=\frac{n+1+\sqrt{n^2+10n-7}}{2(n-1)}. 
\end{equation}
We note that $p_0(n)$ is monotonously decreasing in $n$.
This conjecture had been verified by many authors with partial results.
All the references on the final result in each part
can be summarized in the following table.
\begin{center}
\begin{tabular}{|c||c|c|c|}
\hline
& $p<p_0(n)$ & $p=p_0(n)$ & $p>p_0(n)$ \\
\hline
\hline
$n=2$ & Glassey \cite{G81a} & Schaeffer \cite{Sc85} & Glassey \cite{G81b}\\
\hline
$n=3$ & John \cite{J79} & Schaeffer \cite{Sc85} & John \cite{J79}\\
\hline
$n\ge4$ & Sideris \cite{Si84} &  
$
\begin{array}{l}
\mbox{Yordanov $\&$ Zhang \cite{YZ06}}\\
\mbox{Zhou \cite{Z07}, indep.}
\end{array}
$
&
$
\begin{array}{l}
\mbox{Georgiev $\&$ Lindblad}\\
\mbox{$\&$ Sogge \cite{GLS97}}\\
\end{array}
$
\\
\hline
\end{tabular} 
\end{center}
\par
In the blow-up case, i.e. $1<p\le p_0(n)$,
we are interested in the estimate of the lifespan $T(\e)$.
From now on, $c$ and $C$ stand for positive constants but independent of $\e$.
When $n=1$, we have the following estimate of the lifespan $T(\e)$ for any $p>1$.
\begin{equation}
\label{lifespan_1d}
\left\{
\begin{array}{cl}
c\e^{-(p-1)/2}\le T(\e)\le C\e^{-(p-1)/2}
& \mbox{if}\quad\d\int_{\R}g(x)dx\neq0,\\
c\e^{-p(p-1)/(p+1)}\le T(\e)\le C\e^{-p(p-1)/(p+1)}
&\mbox{if} \quad\d\int_{\R}g(x)dx=0.
\end{array}
\right.
\end{equation}
This result has been obtained by Zhou \cite{Z92_one}.
Moreover, Lindblad \cite{L90} has obtained more precise result for $p=2$,
\begin{equation}
\label{lifespan_1d_lim}
\left\{
\begin{array}{ll}
\d \exists \lim_{\e\rightarrow+0}\e^{1/2}T(\e)>0
&\mbox{for}\quad\d\int_{\R}g(x)dx\neq0,\\
\d \exists \lim_{\e\rightarrow+0}\e^{2/3}T(\e)>0
&\mbox{for}\quad\d\int_{\R}g(x)dx=0.
\end{array}
\right.
\end{equation}
Similarly to this, Lindblad \cite{L90} has also obtained the following result
for $(n,p)=(2,2)$.
\begin{equation}
\label{lifespan_2d_lim}
\left\{
\begin{array}{ll}
\d \exists \lim_{\e\rightarrow+0}a(\e)^{-1}T(\e)>0
&\mbox{for}\quad\d\int_{\R^2}g(x)dx\neq0\\
\d \exists \lim_{\e\rightarrow+0}\e T(\e)>0
&\mbox{for}\quad\d\int_{\R^2}g(x)dx=0,
\end{array}
\right.
\end{equation}
where $a=a(\e)$ is a number satisfying 
\begin{equation}
\label{a}
a^2\e^2\log(1+a)=1.
\end{equation}
\par
When $1<p<p_0(n)\ (n\ge3)$ or $2<p<p_0(2)\ (n=2)$,
we have the following conjecture.
\begin{equation}
\label{lifespan_high-d}
c\e^{-2p(p-1)/\gamma(p,n)}\le T(\e)\le C\e^{-2p(p-1)/\gamma(p,n)},
\end{equation}
where $\gamma(p,n)$ is defined by (\ref{gamma}).
We note that (\ref{lifespan_high-d}) coincides
with the second line in (\ref{lifespan_1d})
if we define $\gamma(p,n)$ by (\ref{gamma}) even for $n=1$.
All the results verifying this conjecture are summarized in the following table.
\begin{center}
\begin{tabular}{|c||c|c|c|}
\hline
& lower bound of $T(\e)$ & upper bound of $T(\e)$\\
\hline
\hline
$n=2$ & Zhou \cite{Z93} & Zhou \cite{Z93}\\
\hline
$n=3$ & Lindblad \cite{L90} & Lindblad \cite{L90} \\
\hline  
$n\ge4$ 
&
Lai $\&$ Zhou \cite{LZ14}
&
$
\begin{array}{l}
\mbox{(rescaling argument}\\
\mbox{of Sideris \cite{Si84})}
\end{array}
$
\\
\hline
\end{tabular} 
\end{center}
We note that, for $n=2,3$,
\[
 \exists \lim_{\e\rightarrow+0}\e^{2p(p-1)/\gamma(p,n)}T(\e)>0
\]
is established in this table.
Moreover, it has been believed that the upper bound in the case where $n\ge 4$
easily follows from the rescaling method applied to the proof in Sideris \cite{Si84}
which proves $T(\e)<\infty$.
Such an argument is actually employed in Georgiev, Takamura and Zhou \cite{GTZ06}
for the analysis on the system and in Zhou and Han \cite{ZH11}
for the analysis on the boundary value problem in the exterior domain.
But it requires more delicate analysis in the following sense.
The conclusion of the original Kato's lemma in \cite{Kato80} is that
the integral in full space of unknown function blows-up in finite time.
The rescaling argument uses this blow-up time as a coefficient in front of
the order of $\e$ with a negative power.
To do this, we have to clarify that such a blow-up time of the rescaled solution does not depend on $\e$
while the initial data of the rescaled solution has some order of $\e$.
Of course, it is true as in Lemma \ref{lem:improvedKato1} in Section \ref{section:ODI} below.
But we do not need the rescaling argument as we see later
if we employ such an improved Kato's lemma.
The purpose of this paper is to show this story.
\par
When $p=p_0(n)$,
we have the following conjecture.
\begin{equation}
\label{lifespan_critical}
\exp\left(c\e^{-p(p-1)}\right)\le T(\e)\le\exp\left(C\e^{-p(p-1)}\right).
\end{equation}
All the results verifying this conjecture are also summarized in the following table.
\begin{center}
\begin{tabular}{|c||c|c|c|}
\hline
& lower bound of $T(\e)$ & upper bound of $T(\e)$\\
\hline
\hline
$n=2$ & Zhou \cite{Z93} & Zhou \cite{Z93}\\
\hline
$n=3$ & Zhou \cite{Z92_three} & Zhou \cite{Z92_three} \\
\hline  
$n\ge4$ 
&
$
\begin{array}{l}
\mbox{Lindblad $\&$ Sogge \cite{LS96}}\\
\mbox{: $n\le 8$ or radially symm. sol.}
\end{array}
$
&
Takamura $\&$ Wakasa \cite{TW11}\\
\hline
\end{tabular} 
\end{center}
Our motivation of this work comes from \cite{TW11},
in which the improved Kato's lemma in the critical case is one of the keys
for the success to finalize this part. 
\par
This paper is organized as follows.
In the next section, we improve Kato's lemma
showing estimates of the existence time of unknown functions
in terms of the holding time of the key inequality.
In Section \ref{section:lifespan}, the improved Kato's lemma is applied to semilinear wave equations.
In Section \ref{section:2d} and \ref{section:1d},
some improvements of the estimates are given in two and one dimensional cases.


\section{Improved Kato's lemma}
\label{section:ODI}
\par
Kato's lemma on ordinary differential inequality in \cite{Kato80} is improved here. 
\begin{lem}
\label{lem:improvedKato1}
Let $p>1, a>0, q>0$ satisfying
\begin{equation}
\label{subcriticalexponent}
M:=\frac{p-1}{2}a-\frac{q}{2}+1>0.
\end{equation}
Assume that $F\in C^2([0,T))$ satisfy
\begin{eqnarray}
& F(t)\ge  At^a &  \mbox{for}\ t\ge T_0,\label{ineq:A}\\
& F''(t)\ge  B(t+R)^{-q}|F(t)|^p & \mbox{for}\ t\ge0,\label{ineq:B}\\
& F(0)\ge0,\ F'(0)>0,\label{ineq:C}&
\end{eqnarray}
where $A,B,R,T_0$ are positive constants.
Then, there exists a positive constant $C_0=C_0(p,a,q,B)$ such that
\begin{equation}
\label{est:T_1}
T<2^{2/M}T_1
\end{equation}
holds provided
\begin{equation}
\label{condi}
T_1:=\max\left\{T_0,\frac{F(0)}{F'(0)},R\right\}\ge C_0 A^{-(p-1)/(2M)}.
\end{equation}
\end{lem}
\begin{rem}
\label{rem:Kato1}
The statement of the original version of this lemma in Kato \cite{Kato80}
is simply that $T<\infty$ without (\ref{ineq:A}),
but there is no restriction on $F(0)$.
The first paper to put an additional assumption (\ref{ineq:A})
and to get a sharp blow-up result of solutions of semilinear wave equations for $n\ge2$
is Glassey \cite{G81a} in two space dimensions.
See also Lemma 4 in Sideris \cite{Si84} for other dimensions.
\end{rem}
\par\noindent
{\bf Proof of Lemma \ref{lem:improvedKato1}.}
Since (\ref{ineq:B}) and (\ref{ineq:C}) imply $F'(t)\ge F'(0)>0$, we have
\begin{equation}
\label{bound:F1}
F(t)\ge F'(0)t+F(0)\ge F(0)\ge0\quad\mbox{for}\ t\ge0.
\end{equation}
Multiplying (\ref{ineq:B}) by $F'(t)>0$, we get
\[
\left(\frac{1}{2}F'(t)^2\right)'\ge B(t+R)^{-q}F(t)^pF'(t)\quad\mbox{for}\ t\ge0.
\] 
It follows from this inequality and (\ref{bound:F1}) that
\[
\begin{array}{ll}
\d\frac{1}{2}F'(t)^2
&\d \ge\frac{1}{2}F'(0)^2+B(t+R)^{-q}\int_0^tF(s)^pF'(s)ds\\
&\d >\frac{B}{p+1}(t+R)^{-q}\left\{F(t)^{p+1}-F(0)^{p+1}\right\}\\
&\d \ge\frac{B}{p+1}(t+R)^{-q}F(t)^p\left\{F(t)-F(0)\right\}
\end{array}
\]
for $t\ge0$.
\par
If $F(0)>0$, then we shall assume that $t\ge F(0)/F'(0)$ which implies $F(t)\ge 2F(0)$.
Hence we have
\[
F'(t)>\sqrt{\frac{B}{p+1}}(t+R)^{-q/2}F(t)^{(p+1)/2}\quad\mbox{for}\ t\ge\frac{F(0)}{F'(0)}.
\]
We note that it is trivial that this inequality also holds for $F(0)=0$.
\par
From now on, we assume that $t\ge T_1$. Then, (\ref{ineq:A}) is available,
so that taking $\delta$ as $0<\delta<(p-1)/2$, we obtain
\[
\begin{array}{ll}
\d\frac{F'(t)}{F(t)^{1+\delta}}
&\d >2^{-q/2}\sqrt{\frac{B}{p+1}}t^{-q/2}F(t)^{(p-1)/2-\delta}\\
&\d >2^{-q/2}\sqrt{\frac{B}{p+1}}A^{(p-1)/2-\delta}t^{M-1-\delta a}
\end{array}
\]
for $t\ge T_1$.
Integrating this inequality on $[T_1,t]$, we get
\[
\frac{1}{\delta}\left(\frac{1}{F(T_1)^\delta}-\frac{1}{F(t)^\delta}\right)
>2^{-q/2}\sqrt{\frac{B}{p+1}}A^{(p-1)/2-\delta}\int_{T_1}^t s^{M-1-\delta a}ds.
\]
Neglecting the second term on the left-hand side
and restricting $\delta$ further as
\[
0<\delta<\min\left(\frac{p-1}{2},\frac{M}{2a}\right),
\]
we have
\[
\frac{1}{\delta F(T_1)^\delta}
>2^{-q/2}\sqrt{\frac{B}{p+1}}A^{(p-1)/2-\delta}\frac{t^{M-\delta a}-T_1^{M-\delta a}}{M-\delta a}
\]
for $t\ge T_1$.
Then, it follows from (\ref{ineq:A}) with $t=T_1$ that
\[
\frac{1}{T_1^{\delta a}}\left(\frac{2^{-q/2}}{M-\delta a}
\sqrt{\frac{B}{p+1}}A^{(p-1)/2}T_1^M+\frac{1}{\delta}\right)
>\frac{2^{-q/2}}{M-\delta a}\sqrt{\frac{B}{p+1}}A^{(p-1)/2}t^{M-\delta a}
\]
for $t\ge T_1$.
If we assume further that
\[
\frac{2^{-q/2}}{M-\delta a}
\sqrt{\frac{B}{p+1}}A^{(p-1)/2}T_1^M\ge\frac{1}{\delta},
\]
namely (\ref{condi}) with
\[
C_0=\left(\frac{2^{-q/2}\delta}{M-\delta a}\sqrt{\frac{B}{p+1}}\right)^{-1/M}>0,
\]
we obtain $2T_1^{M-\delta a}>t^{M-\delta a}$, that is
\[
2^{2/M}T_1>2^{1/(M-\delta a)}T_1>t.
\]
Therefore the proof of the lemma is now completed.
\hfill$\Box$ 
\vskip10pt
We have one more lemma with a different assumption on the initial data.
\begin{lem}
\label{lem:improvedKato2}
Assume that (\ref{ineq:C}) is replaced by
\begin{equation}
\label{ineq:C'}
F(0)>0,\ F'(0)=0
\end{equation}
and additionally that there is a time $t_0>0$ such that  
\begin{equation}
\label{condition_F'}
F(t_0)\ge 2F(0).
\end{equation}
Then, the conclusion of Lemma \ref{lem:improvedKato1} is changed to that
there exists a positive constant $C_0=C_0(p,a,q,B)$ such that
\begin{equation}
\label{est:T_2}
T<2^{2/M}T_2
\end{equation}
holds provided
\begin{equation}
\label{condi2}
T_2:=\max\left\{T_0,t_0,R\right\}\ge C_0 A^{-(p-1)/(2M)}.
\end{equation}
\end{lem}
\begin{rem}
\label{rem:Kato2}
The statement of the original version of this lemma in Kato \cite{Kato80}
is simply that $T<\infty$ without (\ref{ineq:A}), but the assumption on the data is only $F(0)\neq0$
and there is no restriction such as (\ref{condition_F'}).
\end{rem}
\par\noindent
{\bf Proof of Lemma \ref{lem:improvedKato2}.}
It follows from (\ref{ineq:B}) with $t=0$ and (\ref{ineq:C'}) that
$F''(0)>0$ which implies $F'(t)>F'(0)=0$ for small $t>0$. 
Hence $F''(t)\ge0$ for $t\ge0$ from (\ref{ineq:B}) yields that
\[
F'(t)>0\quad\mbox{and}\quad F(t)>F(0)>0\quad\mbox{for}\ t>0.
\]
Hence, similarly to the proof of Lemma \ref{lem:improvedKato1},
we obtain that
\[
\frac{1}{2}F'(t)^2\ge\frac{B}{p+1}(t+R)^{-q}F(t)^p\left\{F(t)-F(0)\right\}\quad\mbox{for}\ t>0.
\]
Then taking into account of (\ref{condition_F'}) and $F(t)\ge F(t_0)$ for $t\ge t_0$,
we have that
\[
F'(t)>\sqrt{\frac{B}{p+1}}(t+R)^{-q/2}F(t)^{(p+1)/2}
\quad\mbox{for}\ t\ge t_0.
\]
Therefore, after assuming $t\ge T_2$, one can readily establish the same argument
as in the proof of Lemma \ref{lem:improvedKato1} in which
$T_1$ is replaced by $T_2$.
\hfill$\Box$ 

\section{Upper bound of the lifespan}
\label{section:lifespan}
In this section, we prove the expected theorem on the upper bound of the lifespan
in high dimensional case.
\begin{thm}
\label{thm:main1}
Let $1<p<p_0(n)$ for $n\ge2$.
Assume that both $f\in H^1(\R^n)$ and $g\in L^2(\R^n)$
are non-negative and have compact support,
and that $g$ does not vanish identically.
Suppose that the problem {\rm(\ref{IVP})} has a solution
$(u,u_t)\in C([0,T(\e)),H^1(\R^n)\times L^2(\R^n))$ with
\begin{equation}
\label{support}
\mbox{\rm supp}(u,u_t)\subset\{(x,t)\in\R^n\times[0,\infty)\ :\ |x|\le t+R\}.
\end{equation}
Then, there exists a positive constant $\e_0
=\e_0(f,g,n,p,R)$ such that $T(\e)$ has to satisfy
\begin{equation}
\label{lifespan1}
T(\e)\le C\e^{-2p(p-1)/\gamma(p,n)}
\end{equation}
for $0<\e\le\e_0$, where $C$ is a positive constant independent of $\e$.
\end{thm}
\begin{rem}
\label{rem:thm1}
In view of (\ref{lifespan_2d_lim}), (\ref{lifespan1}) is not optimal for $n=2$ and $1<p\le2$.
In Section \ref{section:2d} and \ref{section:1d}, by improving estimate  (\ref{step0}),
optimal estimates for the upper bound for the lifespan are obtained for $n=2$ and $n=1$.
\end{rem}
\par\noindent
{\bf Proof of Theorem \ref{thm:main1}.}
In order to employ Lemma \ref{lem:improvedKato1}, let us set
\[
F(t)=\int_{\R^n}u(x,t)dx.
\]
Then, (\ref{ineq:C}) immediately follows from the assumption of the theorem as
\begin{equation}
\label{initial}
F(0)=\e\int_{\R^n}f(x)dx\ge0,\quad F'(0)=\e\int_{\R^n}g(x)dx>0.
\end{equation}
\par
For (\ref{ineq:B}), we shall employ the same argument as (13)-(15) in Sideris \cite{Si84}.
Then it follows that
\[
F''(t)=\int_{\R^n}|u(x,t)|^pdx
\ge\frac{|F(t)|^p}{\d\left(\int_{|x|\le t+R}1dx\right)^{p-1}}.
\]
This inequality means that (\ref{ineq:B}) is available with
\begin{equation}
\label{Bq}
B=\left(\int_{|x|\le1}dx\right)^{1-p}>0,\quad q=n(p-1)>0.
\end{equation}
\par
For the key inequality (\ref{ineq:A}), we employ the following proposition.
\begin{prop}
\label{prop:step0}
Suppose that the assumption in Theorem \ref{thm:main1} is fulfilled.
Then, there exists a positive constant $C_1=C_1(f,g,n,p,R)$
such that $F(t)=\d\int_{\R^n}\!\!\!u(x,t)dx$ satisfies
\begin{equation}
\label{step0}
F''(t)\ge C_1\e^pt^{(n-1)(1-p/2)}\quad\mbox{for}\ t\ge R.
\end{equation}
This estimate is valid also for the case where $f(x)\ge0(\not\equiv0)$ and $g(x)\equiv0$.
\end{prop}
\begin{rem}
This is a slightly modified estimate of (2.5') in Yordanov and Zhang \cite{YZ06}
in which $C_1=0$ when $f(x)\equiv0$.
\end{rem}
\par\noindent
{\bf Proof of Proposition \ref{prop:step0}.}
It follows from Lemma 2.2 in \cite{YZ06} that
\[
F_1(t)\ge\frac{1}{2}(1-e^{-2R})\int_{\R^n}\{\e f(x)+\e g(x)\}\phi_1(x)dx
\quad\mbox{for}\ t\ge R,
\]
where
\[
\phi_1(x)=\int_{S^{n-1}}e^{x\cdot\omega}d\omega
\quad\mbox{and}\quad
F_1(t)=\int_{\R^n}u(x,t)\phi_1(x)e^{-t}dx.
\]
Combining this estimate and (2.4) in \cite{YZ06}, we immediately obtain  (\ref{step0}).
We note that there is no restriction on $n$ in the argument.
\hfill$\Box$
\vskip10pt
\par
Let us continue to prove the theorem.
Integrating (\ref{step0}) in $[R,t]$, we have
\[
F'(t)\ge\frac{C_1}{n-(n-1)p/2}\e^p\left(t^{n-(n-1)p/2}-R^{n-(n-1)p/2}\right)+F'(R)
\]
for $t\ge R$ because of
\[
1<p<p_0(n)<\frac{2n}{n-1}\quad\mbox{for}\ n\ge2.
\]
Note that $F'(R)>0$ follows from $F''(t)\ge0$ for $t\ge0$ and $F'(0)>0$.
Hence we obtain that
\begin{equation}
\label{step1}
F'(t)>\frac{C_1\left(1-2^{-n+(n-1)p/2}\right)}{n-(n-1)p/2}\e^pt^{n-(n-1)p/2}
\quad\mbox{for}\ t\ge2R.
\end{equation}
Integrating this inequality in $[2R,t]$ together with $F(0)\ge0$, we get
\begin{equation}
\label{condition_F}
F(t)>C_2\e^pt^{n+1-(n-1)p/2}
\quad\mbox{for}\ t\ge4R,
\end{equation}
where
\[
C_2=\frac{C_1\left(1-2^{-n+(n-1)p/2}\right)\left(1-2^{-n-1+(n-1)p/2}\right)}
{\left\{n-(n-1)p/2\right\}\left\{n+1-(n-1)/2\right\}}>0.
\]
\par
We are now in a position to apply our result here to Lemma \ref{lem:improvedKato1}
with special choices on all positive constants except for $T_0$ as
\[
A=C_2\e^p,\ B=\left(\int_{|x|\le1}dx\right)^{1-p},
\ a=n+1-\frac{n-1}{2}p,\ q=n(p-1)
\]
which imply that (\ref{subcriticalexponent}) yields
\[
M=\frac{p-1}{2}a-\frac{q}{2}+1=\frac{\gamma(p,n)}{4}>0.
\]
If we set
\begin{equation}
\label{condition_T_0}
T_0=C_0A^{-(p-1)/(2M)}=C_0C_2^{-2(p-1)/\gamma(p,n)}\e^{-2p(p-1)/\gamma(p,n)},
\end{equation}
we then find that there is an $\e_0=\e_0(f,g,n,p,R)>0$ such that
\[
T_0\ge\max\left\{\frac{F(0)}{F'(0)},4R\right\}
\quad\mbox{for}\ 0<\e\le\e_0
\] 
because $F(0)/F'(0)$ does not depend on $\e$.
This means that $T_1=T_0$ in (\ref{condi}).
Therefore the conclusion of Lemma \ref{lem:improvedKato1} implies that
the maximal existence time $T$ of $F(t)$ has to satisfy
\[
T(\e)\le T\le C_3\e^{-2p(p-1)/\gamma(p,n)}
\quad\mbox{for}\ 0<\e\le\e_0,
\]
where
\[
C_3=2^{8/\gamma(p,n)}C_0C_2^{-2(p-1)/\gamma(p,n)}>0.
\]
The proof of the theorem is now completed.
\hfill$\Box$
\begin{thm}
\label{thm:main2}
Let $1<p<p_0(n)$ for $n\ge2$.
Assume that $f\in H^1(\R^n)$ is non-negative, does not vanish identically and $g\equiv0$. 
Suppose that the problem {\rm(\ref{IVP})} has a solution
$(u,u_t)\in C([0,T(\e)),H^1(\R^n)\times L^2(\R^n))$ with (\ref{support}).
Then, there exists a positive constant $\e_1=\e_1(f,n,p,R)$ such that $T(\e)$ has to satisfy
\begin{equation}
\label{lifespan2}
T(\e)\le C\e^{-2p(p-1)/\gamma(p,n)}
\end{equation}
for $0<\e\le\e_1$, where $C$ is a positive constant independent of $\e$.
\end{thm}
\begin{rem}
We have no new estimate for high dimensional case in this theorem,
but it especially covers the optimality on $n=2$ and $p=2$.
See (\ref{lifespan_2d_lim}).
\end{rem}
\par\noindent
{\bf Proof of Theorem \ref{thm:main2}.}
Note that Proposition \ref{prop:step0} is also available in this case.
Hence we have (\ref{condition_F}) again.
Assume that
\[
C_2\e^pt_0^{n+1-(n-1)p/2}
=2F(0)=2\|f\|_{L^1(\R^n)}\e
\quad\mbox{and}\quad t_0\ge4R.
\]
Then (\ref{condition_F'}) in Lemma \ref{lem:improvedKato2} is fulfilled with
\begin{equation}
\label{t_0}
t_0=C_4\e^{-(p-1)/(n-(n-1)p/2)}
\end{equation}
if $t_0\ge4R$, where
\[
C_4=\left\{
\frac{2\|f\|_{L^1(\R^n)}}{C_2}
\right\}^{\left\{n+1-(n-1)p/2\right\}^{-1}}>0.
\]
\par
Since
\[
\frac{1}{n+1-(n-1)p/2}\le\frac{2p}{\gamma(p,n)}
\]
is equivalent to the trivial condition:
\[
p\ge\frac{2}{n+1},
\]
we obtain that there is an $\e_1=\e_1(f,n,p,R)>0$ such that
\[
T_0\ge t_0\quad\mbox{and}\quad t_0\ge4R\quad\mbox{for}\ 0<\e\le\e_1
\]
with the same choice of $T_0$ as in (\ref{condition_T_0}).
This means $T_2=T_0$ in (\ref{condi2}),
so that the same conclusion as in Theorem \ref{thm:main2} holds.
\hfill$\Box$
\section{Note on 2 dimensional case}
\label{section:2d}
The optimal estimate in $n=2$ and $p=2$ in the case where $g(x)\ge0(\not\equiv0)$
is obtained by improving (\ref{step0}) as announced in Remark \ref{rem:thm1}.
 
\begin{thm}
\label{thm:2d}
Let $n=2$, $1<p\le2$ and $f\equiv0$.
Assume that $g\in C^1(\R^2)$ is non-negative, does not vanish identically, and 
has compact support as supp $g\subset\{x\in\R^2\ :\ |x|\le R\}$.
Suppose that the integral equation associated with the problem {\rm(\ref{IVP})} has a solution
$u\in C^1(\R^2\times[0,T(\e)))$ with supp $u(x,t)\subset\{(x,t)\in\R^2\times[0,\infty)\ :\ |x|\le t+R\}$.
Then, there exists a positive constant $\e_2
=\e_2(g,p,R)$ such that $T(\e)$ has to satisfy
\begin{equation}
\label{lifespan_2d}
T(\e)\le
\left\{
\begin{array}{ll}
C a(\e) & \mbox{when}\ p=2,\\
C\e^{-(p-1)/(3-p)} & \mbox{when}\ 1<p<2
\end{array}
\right.
\end{equation}
for $0<\e\le\e_2$, where $a(\e)$ is defined in (\ref{a})
and  $C$ is a positive constant independent of $\e$.
\end{thm}
\begin{rem}
(\ref{lifespan_2d}) is a weaker result than (\ref{lifespan_2d_lim})
in the sense that the constant $C$ is not of the optimal choice
and that the assumption on the data is stronger than Lindblad \cite{L90}.
Moreover, we note that the optimality of (\ref{lifespan_2d})
for $1<p<2$ is still open, but may true.
Because
\[
\frac{1}{3-p}<\frac{2p}{\gamma(p,2)}
\]
is equivalent to $1<p<2$.
But it is out of the main purpose of this paper.
So we shall have another opportunity to study this part. 
\end{rem}
\par\noindent
{\bf Proof of Theorem \ref{thm:2d}.}
Under the assumption of the theorem,
it is well-known that our integral equation is of the form,
\begin{equation}
\label{rep_2d}
\begin{array}{ll}
u(x,t)=
&\d\frac{\e}{2\pi}\int_{|x-y|\le t}\frac{g(y)}{\sqrt{t^2-|x-y|^2}}dy\\
&\d+\frac{1}{2\pi}\int_0^td\tau
\int_{|x-y|\le t-\tau}\frac{|u(y,\tau)|^p}{\sqrt{(t-\tau)^2-|x-y|^2}}dy.
\end{array}
\end{equation}
Note that $|y|\le R$ and $|x|\le t+R$ due to the support property
in the first term on the right-hand side.
Then, neglecting the second term on the right-hand side
and making use of inequalities;
\[
\begin{array}{l}
t-|x-y|\le t-\left||x|-|y|\right|\le t-|x|+R\quad\mbox{for}\ |x|\ge R,\\
t+|x-y|\le t+|x|+R\le2(t+R),
\end{array}
\]
we obtain that
\[
u(x,t)\ge\frac{\e}{2\sqrt{2}\pi\sqrt{t+R}\sqrt{t-|x|+R}}\int_{|x-y|\le t}g(y)dy
\quad\mbox{for}\ |x|\ge R.
\]
If we assume $|x|+R\le t$ which implies $|x-y|\le t$ for $|y|\le R$,
we get
\[
\int_{|x-y|\le t}g(y)dy=\|g\|_{L^1(\R^2)}.
\]
Therefore we have that
\begin{equation}
\label{est_u_2d}
 u(x,t)\ge\frac{\|g\|_{L^1(\R^2)}}{2\sqrt{2}\pi\sqrt{t+R}\sqrt{t-|x|+R}}\e
\quad\mbox{for}\ R\le|x|\le t-R.
\end{equation}
\par
It follows from the same argument as in high dimensional case that $\d F(t)=\int_{\R^2}u(x,t)dx$ satisfies
\[
F''(t)=\int_{\R^2}|u(x,t)|^pdx
\ge\int_{R\le|x|\le t-R}|u(x,t)|^pdx\quad\mbox{for}\ t\ge2R.
\]
Plugging (\ref{est_u_2d}) into the right-hand side of this inequality, we have that
\[ 
F''(t)\ge\left(\frac{\|g\|_{L^1(\R^2)}}{2\sqrt{2}\pi\sqrt{t+R}}\e\right)^p
\int_{R\le|x|\le t-R}\frac{1}{(t-|x|+R)^{p/2}}dx,
\]
that is
\begin{equation}
\label{est_F''_2d}
F''(t)\ge\frac{2\pi\|g\|^p_{L^1(\R^2)}}{(2\sqrt{2}\pi)^p(t+R)^{p/2}}\e^p
\int_R^{t-R}\frac{r}{(t-r+R)^{p/2}}dr
\quad\mbox{for}\ t\ge2R.
\end{equation}
\par\noindent
{\bf Case of $\v{p=2}$.} Making use of integration by parts, we have that
\[
\int_R^{t-R}\frac{r}{t-r+R}dr
=R\log t-(t-R)\log2R+\int_R^{t-R}\log(t-r+R)dr.
\]
Without loss of the generality, we may assume that
\begin{equation}
\label{R}
R\ge1.
\end{equation}
Note that $(t-R)/2\ge R$ is equivalent to $t\ge3R$.
Hence, diminishing the domain of the integral to $[R,(t-R)/2]$, we get that
\[
 \int_R^{t-R}\log(t-r+R)dr
\ge\left(\frac{t-R}{2}-R\right)\log\frac{t+3R}{2}
\quad\mbox{for}\ t\ge3R.
\]
Therefore we obtain that
\[
\begin{array}{ll}
\d\int_R^{t-R}\frac{r}{t-r+R}dr\ge
&\d R\left(\log t-\log\frac{t+3R}{2}\right)\\
&\d +(t-R)\left(\frac{1}{2}\log\frac{t+3R}{2}-\log2R\right)
\end{array}
\mbox{for}\ t\ge3R.
\]
Then it follows from this inequality,
\[
\log\frac{t+3R}{2}\le\log t\quad\mbox{for}\ t\ge3R
\]
and
\[
\frac{1}{4}\log\frac{t+3R}{2}\ge\log2R
\quad\mbox{for}\ t\ge32R^4-3R
\]
that
\[
\int_R^{t-R}\frac{r}{t-r+R}dr\ge\frac{1}{6}t\log\frac{t+3R}{2}
\quad\mbox{for}\ t\ge32R^4-3R
\]
because (\ref{R}) implies $32R^4-3R\ge3R$.
Combining this estimate with (\ref{est_F''_2d}), one can obtain that
\[
F''(t)\ge\frac{\|g\|^2_{L^1(\R^2)}}{48\pi}\e^2\log\frac{t}{2}
\quad\mbox{for}\ t\ge 32R^4.
\]
This is a better estimate than (\ref{step0}) as the extra factor $\log(t/2)$ appears,
which leads to the optimal order of the lifespan in this case as follows.
\par
Integrating the inequality above in $[32R^4,t]$ and making use of $F'(0)=\|g\|_{L^1(\R^2)}\e>0$,
we get
\[
F'(t)>\frac{\|g\|^2_{L^1(\R^2)}}{48\pi}\e^2\int_{t/2}^t\log\frac{s}{2}ds
\ge\frac{\|g\|^2_{L^1(\R^2)}}{96\pi}\e^2t\log\frac{t}{4}
\quad\mbox{for}\ t\ge 64R^4.
\]
The same procedure with $F(0)=0$ in $[64R^2,t]$ gives us that
\begin{equation}
\label{est:F}
F(t)>\frac{\|g\|^2_{L^1(\R^2)}}{96\pi}\e^2\int_{t/2}^ts\log\frac{s}{4}ds
\ge\frac{\|g\|^2_{L^1(\R^2)}}{256\pi}\e^2t^2\log\frac{t}{8}
\quad\mbox{for}\ t\ge 128R^4.
\end{equation}
We are now in a position to apply our situation to Lemma \ref{lem:improvedKato1} with
\[
a=p=q=2,B=\pi^{-1}
\]
which imply $M=\gamma(2,2)/4=1$.
If $T_0$ satisfies
\begin{equation}
\label{est:T_0}
T_0\ge\max\left\{\frac{F(0)}{F'(0)},128R^4\right\},
\end{equation}
then we have $T_1=T_0$ which yields that
 (\ref{condi}) is equivalent to
\begin{equation}
\label{condi:C_0}
A\ge C_0^2T_0^{-2}.
\end{equation}
Moreover,  in view of (\ref{est:F}), (\ref{ineq:A}) follows from
\begin{equation}
\label{condi:A}
\frac{\|g\|^2_{L^1(\R^2)}}{256\pi}\e^2\log\frac{T_0}{8}\ge A.
\end{equation}
Hence the requirement to $T_0$ is
\begin{equation}
\label{condi:T_0}
\frac{\|g\|^2_{L^1(\R^2)}}{256\pi}T_0^2\e^2\log\frac{T_0}{8}\ge C_0^2.
\end{equation}
\par
For this possibility on the choice of $T_0$, one can set
\[
\frac{T_0}{16}=a\left(\e\right)\quad\mbox{when}\quad \frac{\|g\|^2_{L^1(\R^2)}}{\pi C_0^2}\ge1,
\]
where $a(\e)$ is defined in (\ref{a}).
Because there exists a positive constant $\e_{21}=\e_{21}(g,R)$ such that
\[
a(\e)\ge\max\left\{\frac{F(0)}{16F'(0)},8R^4\right\}
\quad\mbox{for}\ 0<\e\le\e_{21}
\]
which yields (\ref{est:T_0})
due to the fact that $a(\e)$ is a monotonously decreasing function of $\e$
and $\d\lim_{\e\rightarrow+0}a(\e)=\infty$.
Then it follows from (\ref{R}) that
\[
\frac{T_0}{8}=2a(\e)\ge8R^4+a(\e)\ge1+a(\e)
\quad\mbox{for}\ 0<\e\le\e_{21}.
\]
Therefore (\ref{est:T_0}) and (\ref{condi:T_0}) are established by setting
\begin{equation}
\label{choice1:T_0}
T_0=16a(\e)\quad\mbox{for}\quad0<\e\le\e_{21}
\end{equation}
because we have
\[
\frac{\|g\|^2_{L^1(\R^2)}}{256\pi C_0^2}T_0^2\e^2\log\frac{T_0}{8}
\ge\left(\frac{T_0}{16}\right)^2\e^2\log\frac{T_0}{8}
\ge a(\e)^2\e^2\log\left(1+a(\e)\right)=1
\]
for $0<\e\le\e_{21}$.
On the counter part, one can set
\[
\frac{\|g\|_{L^1(\R^2)}}{\sqrt{\pi}C_0}\cdot
\frac{T_0}{16}=a\left(\e\right)\quad\mbox{when}
\quad\frac{\|g\|^2_{L^1(\R^2)}}{\pi C_0^2}\le1.
\]
Because there exists a positive constant $\e_{22}=\e_{22}(g,R)$ such that
\[
a(\e)\ge\frac{\|g\|_{L^1(\R^2)}}{\sqrt{\pi}C_0}\max\left\{\frac{F(0)}{16F'(0)},8R^4\right\}
\quad\mbox{for}\ 0<\e\le\e_{22}
\]
which yields (\ref{est:T_0}).
Then it follows from (\ref{R}) that
\[
\frac{T_0}{8}=\frac{2\sqrt{\pi}C_0}{\|g\|_{L^1(\R^2)}}a(\e)
\ge8R^4+\frac{\sqrt{\pi}C_0}{\|g\|_{L^1(\R^2)}}a(\e)\ge 1+a(\e)\quad\mbox{for}\ 0<\e\le\e_{22}.
\]
Therefore (\ref{est:T_0}) and (\ref{condi:T_0}) are established by setting
\begin{equation}
\label{choice2:T_0}
T_0=\frac{16\sqrt{\pi}C_0}{\|g\|_{L^1(\R^2)}}a(\e)\quad\mbox{for}\quad0<\e\le\e_{22}
\end{equation}
because we have
\[
\frac{\|g\|^2_{L^1(\R^2)}}{256\pi C_0^2}T_0^2\e^2\log\frac{T_0}{8}
\ge a(\e)^2\e^2\log\left(1+a(\e)\right)=1
\]
for $0<\e\le\e_{22}$.
Summing up (\ref{choice1:T_0}) and (\ref{choice2:T_0}),
we can set
\[
T_0=16\max\left\{1,\frac{\sqrt{\pi}C_0}{\|g\|_{L^1(\R^2)}}\right\}a(\e)
\ \mbox{with}\ 0<\e\le\e_2
\quad\mbox{and}\ A=C_0^2T_0^{-2}
\]
in Lemma \ref{lem:improvedKato1}, where $\e_2=\min\{\e_{21},\e_{22}\}$.
Then the conclusion with $2^{2/M}T_1=4T_0$ yields that
\[
T(\e)\le64\max\left\{1,\frac{\sqrt{\pi}C_0}{\|g\|_{L^1(\R^2)}}\right\}a(\e)
\quad\mbox{for}\ 0<\e\le\e_2.
\]
\par\noindent
{\bf Case of $\v{1<p<2}$.} Turning back to (\ref{est_F''_2d}),
we find that
\[
\int_R^{t-R}\frac{r}{(t-r+R)^{p/2}}dr
\ge\frac{1}{2t^{p/2}}\left\{(t-R)^2-R^2\right\}.
\]
Hence it follows that
\[
\int_R^{t-R}\frac{r}{(t-r+R)^{p/2}}dr
\ge\frac{1}{6}t^{2-p/2}\quad\mbox{for}\ t\ge3R
\]
which yields
\[
F''(t)\ge\frac{\|g\|^p_{L^1(\R^2)}}{3\cdot2^{2p}\pi^{p-1}}\e^pt^{2-p}
\quad\mbox{for}\ t\ge3R.
\]
We note that this is a better estimate than (\ref{step0}) because
\[
1-\frac{p}{2}<2-p
\]
is equivalent to $p<2$.
Integrating this inequality in $[3R,t]$
and making use of $F'(0)=\|g\|_{L^1(\R^2)}\e>0$, we get
\[
F'(t)>\frac{\|g\|^p_{L^1(\R^2)}(1-(3/4)^{3-p})}{3(3-p)2^{2p}\pi^{p-1}}\e^pt^{3-p}
\quad\mbox{for}\ t\ge4R.
\]
Hence we obtain that
\[
F(t)>C_5\e^pt^{4-p}
\quad\mbox{for}\ t\ge5R,
\]
where
\[
C_5=\frac{\|g\|^p_{L^1(\R^2)}(1-(3/4)^{3-p})(1-(4/5)^{4-p})}
{3(3-p)(4-p)2^{2p}\pi^{p-1}}>0.
\]
\par
We are in a position to apply our situation to Lemma \ref{lem:improvedKato1} with
\[
A=C_5\e^p,\ B=\pi^{1-p},\ a=4-p,\ q=2(p-1)
\]
which imply
\[
M=\frac{p-1}{2}(4-p)-\frac{2(p-1)}{2}+1=\frac{p(3-p)}{2}>0.
\] 
Therefore the theorem follows from setting
\[
T_0=C_0A^{-(p-1)/(2M)}=C_0C_5^{-(p-1)/(p(3-p))}\e^{-(p-1)/(3-p)}
\]
as in the proof of Theorem \ref{thm:main1}.
\hfill$\Box$


\section{Note on one dimensional case}
\label{section:1d}
The optimal estimate in $n=1$ in the case where $g(x)\ge0(\not\equiv0)$
is also obtained by a better estimate than (\ref{step0}) as announced in Remark \ref{rem:thm1}.
 
\begin{thm}
\label{thm:1d}
Let $p>1$ for $n=1$.
Assume that both $f\in C^2(\R)$ and $g\in C^1(\R)$
have compact support as {\rm supp} $(f,g)\subset\{x\in\R\ :\ |x|\le R\}$.
Suppose that the problem {\rm(\ref{IVP})} has a solution
$u\in C^2(\R\times[0,T(\e)))$
Then, there exists a positive constant $\e_3
=\e_3(f,g,p,R)$ such that $T(\e)$ has to satisfy
\begin{equation}
\label{lifespan_1d_u}
\begin{array}{ll}
\d T(\e)\le C \e^{-(p-1)/2} & \d\mbox{if}\ \int_{\R}g(x)dx>0,\\
\d T(\e)\le C\e^{-p(p-1)/(p+1)} &\d\mbox{if}\ f\ge0(\not\equiv0),\ g\equiv0
\end{array}
\end{equation}
for $0<\e\le\e_3$, where $C$ is a positive constant independent of $\e$.
\end{thm}
\begin{rem}
The assumption of this theorem is stronger than Zhou \cite{Z92_one}.
But for the sake of completeness of applications of our lemma, we shall prove it here.
\end{rem}
\par\noindent
{\bf Proof of Theorem \ref{thm:1d}.}
First, we note that the assumption on the support of the initial data implies that
the solution $u\in C^2(\R\times[0,T(\e))$
of (\ref{IVP}) has to satisfy
\[
\mbox{supp}\ u(x,t)\subset\{x\in\R\ :\ |x|\le t+R\}.
\]
Hence, integrating the equation in $\R$, we have that $\d F(t)=\int_{\R}u(x,t)dx$ satisfies
\begin{equation}
\label{F''_1d}
F''(t)=\int_{\R}|u(x,t)|^pdx=\int_{-(t+R)}^{t+R}|u(x,t)|^pdx
\end{equation}
which yields
\begin{equation}
\label{ineq:B_1d}
F''(t)\ge2^{1-p}(t+R)^{1-p}|F(t)|^p
\quad\mbox{for}\ t\ge0.
\end{equation}
It is well-known that $u$ has a representation formula of the form,
\begin{equation}
\label{rep_1d}
\begin{array}{ll}
u(x,t)=
&\d\frac{f(x+t)+f(x-t)}{2}\e+\frac{\e}{2}\int_{x-t}^{x+t}g(\xi)d\xi\\
&\d+\frac{1}{2}\int_0^td\tau\int_{x-t+\tau}^{x+t-\tau}|u(\xi,\tau)|^pd\xi.
\end{array}
\end{equation}
\vskip10pt
\par\noindent
{\bf Case of $\v{\d\int_{\R}g(x)dx>0}$.}
In view of (\ref{rep_1d}),  we find that the support condition implies
\[
u(x,t)\ge G\e\quad\mbox{for}\ x+t\ge R\ \mbox{and}\ x-t\le-R,
\]
where
\[
G=\frac{1}{2}\int_{\R}g(x)dx>0.
\]
Plugging this estimate into (\ref{F''_1d}), we have that
\[
F''(t)\ge\left(G\e\right)^p\int_0^{t-R}dx=G^p\e^p(t-R)
\quad\mbox{for}\ t\ge R.
\]
This is a better estimate than (\ref{step0}).
Integrating this inequality twice in $[R,t]$ and making use of
\[
F'(R)\ge F'(0)>0,\quad F(R)>F(0)\ge0,
\]
we obtain that
\[
F(t)>\frac{G^p}{6}\e^p(t-R)^3
\quad\mbox{for}\ t\ge R
\]
which implies that
\[
F(t)>\frac{G^p}{48}\e^pt^3
\quad\mbox{for}\ t\ge2R.
\]
\par
In view of this estimate and (\ref{ineq:B_1d}),
we are in a position to apply Lemma \ref{lem:improvedKato1} with
\[
A=\frac{G^p}{48}\e^p,B=2^{1-p},a=3,q=p-1.
\]
Noticing that
\[
M=\frac{p-1}{2}3-\frac{p-1}{2}+1=p,
\]
one can set
\[
T_0=C_0A^{-(p-1)/(2M)}=C_0\frac{G^p}{48}\e^{-(p-1)/(2p)}.
\]
Therefore the first conclusion of the theorem follows.
\vskip10pt
\par\noindent
{\bf Case of $\v{f\ge0(\not\equiv0),g\equiv0}$.}
In view of (\ref{rep_1d}),  we find that the support condition implies
\[
u(x,t)\ge\frac{f(x-t)}{2}\e\quad\mbox{for}\ x+t\ge R\ \mbox{and}\ -R\le x-t\le R.
\]
Plugging this estimate into (\ref{F''_1d}), we have that
\[
F''(t)\ge\frac{\e^p}{2^p}\int_{t-R}^{t+R}f(x-t)^pdx
=\frac{\|f\|^p_{L^p(\R)}}{2^p}\e^p
\quad\mbox{for}\ t\ge R.
\]
This is the same estimate as (\ref{step0}).
Therefore the second conclusion of the theorem immediately follows
in the same manner as in the proof of Theorem \ref{thm:main2} because of $\gamma(p,1)=2+2p$.
\hfill$\Box$
\vskip10pt
\par\noindent
{\bf Acknowledgment.} This work is partially supported by
the Grant-in-Aid for Scientific Research (C) (No. 24540183),
Japan Society for the Promotion of Science.
The author is grateful to members of a private seminar
held at Future University Hakodate on 24-25 Jan. 2015,
Professors A.Hoshiga, K.Yokoyama, Y.Kurokawa and Dr. K.Wakasa,
for various discussions which lead to the optimal condition (\ref{condition_F'}) in assumptions of Lemma \ref{lem:improvedKato2}.

\bibliographystyle{plain}

\end{document}